\documentclass[10pt]{amsart}
\usepackage{amssymb, mathrsfs,amsmath}
\usepackage{latexsym}
\usepackage{amsfonts}
\usepackage{shadow}

\newtheorem{Pa}{Paper}[section]
\newtheorem{Tm}[Pa]{{\bf Theorem}}
\newtheorem{La}[Pa]{{\bf Lemma}}

\newtheorem{Ob}[Pa]{{\bf Observation}}

\newtheorem{Pn}[Pa]{{\bf Proposition}}

\newtheorem{Ex}[Pa]{{\bf Example}}

\makeatletter
\def\Ddots{\mathinner{\mkern1mu\raise\p@
\vbox{\kern7\p@\hbox{.}}\mkern2mu
\raise4\p@\hbox{.}\mkern2mu\raise7\p@\hbox{.}\mkern1mu}}
\makeatother
\def\C{\mathbb C}

\date{}
\author[D. Alpay]{Daniel Alpay}
\author[P. Jorgensen]{Palle Jorgensen}
\address{(DA) Department of Mathematics
\newline
Ben Gurion University of the Negev \newline P.O.B. 653,
\newline
Be'er Sheva 84105, \newline ISRAEL} \email{dany@math.bgu.ac.il}
\address{(PJ)
Department of Mathematics\newline 14 MLH \newline The University
of Iowa, Iowa City,\newline IA 52242-1419 USA}
\email{jorgen@math.uiowa.edu}
\author[I. Lewkowicz]{Izchak Lewkowicz}
\address{(IL) Department of Electrical Engineering
\newline
Ben Gurion University of the Negev \newline P.O.B. 653,
\newline
Be'er Sheva 84105, \newline ISRAEL}
\email{izchak@ee.bgu.ac.il}
\thanks{This research is partially supported by the BSF
grant no. 2010117.}
\thanks{D. Alpay thanks the
Earl Katz family for endowing the chair which supported his
research.}
\title
[
(co)-isometric rational Functions
]
{
characterizations of rectangular\vskip 0.3cm
(para)-unitary rational Functions}

\begin{document}
\begin{abstract}
We here present three characterizations of not necessarily causal,
rational functions which are (co)-isometric on the unit circle:\\
(i) Through the realization matrix of Schur stable systems.\\
(ii) The Blaschke-Potapov product, which
is then employed to
introduce an easy-to-use description of all these functions
with dimensions and McMillan degree as parameters.\\
(iii) Through the (not necessarily reducible) Matrix Fraction
Description (MFD).

In cases (ii) and (iii) the poles of the rational functions
involved may be anywhere in the complex plane, but
the unit circle (including both zero and infinity).

A special attention is devoted to exploring the gap between the
square and rectangular cases.
\end{abstract}
\keywords{isometry, coisometry, lossless, all-pass, realization, 
gramians, Matrix Fraction Description, Blaschke-Potapov product.\\
{\em AMS 2010 subject classification index}:
20H05, 26C15, 47A48, 47A56, 51F25, 93B20,
94A05, 94A08, 94A11, 94A12
}\maketitle

\section{Introduction}
\setcounter{equation}{0}
\label{sec:intoduction}

This work is on the crossroads of Operator and Systems theory
from the mathematical side and Control, Signal Processing and
Communications theory from the engineering side. It addresses
problems or employs tools from all these areas. Thus, it is
meant to serve as a bridge between the corresponding communities.
We start by formally laying out the set-up.

\subsection{(Para)-Unitary
symmetry}
\label{subsec:UnitarySym}

Let $F(z)$ be $p\times m$-valued rational functions with poles
everywhere in the complex plane $\C$ (including infinity),
i.e. it can be written as
\begin{equation}\label{eq:rational}
F(z)=C(zI-A)^{-1}B+D+\sum\limits_{j=1}^kz^jE_j~,
\quad\quad\quad\quad k\geq 0,
\end{equation}
where the constant matrices $A$, $B$, $C$, and $D,~E_1~,~\ldots~,~E_k$
are of dimensions $n\times n$, $m\times n$, $p\times n~$ and
$~p\times m$, respectively. Whenever, $k\geq 1$, in system theory
``dialect" $F(z)$ is said to have ~{\em poles at infinity} while in
engineering ``dialect" $F(z)$ is called an ~{\em improper}~ rational
function. Furthermore, $F(z)$ may be viewed as the (two sided)
$Z$-transform of an impulse response $\Phi(t)$, with $t$ an integral
variable. In particular, $k\geq 1$ means that $\Phi(t)\not\equiv 0$
for $0>t$. Thus engineers call it ~{\em non-causal}.
\vskip 0.2cm

Let ${\mathbb T}$ be the unit circle,
\[
{\mathbb T}:=\{z\in\C~:~|z|=1~\}.
\]
In this work we focus on $~\mathcal{U}$, the subclass of
$p\times m$-valued rational functions in \eqref{eq:rational} having
unitary symmetry on the unit circle, i.e.
\begin{equation}\label{eq:DefU}
\mathcal{U}:=\left\{~F(z)~:~\left\{
\begin{matrix}\left(F(z)\right)^*F(z)\equiv I_m&p\geq m&{\rm isometry}\\~\\
F(z)\left(F(z)\right)^*\equiv I_p&m\geq p&{\rm coisometry}\end{matrix}\right.
\quad\forall z\in{\mathbb T}\right\}.
\end{equation}
In signal processing ``dialect"~ {\em unitary}~ is reserved to
constant matrices while ~{\em para-unitary}~ means matrix-valued
functions with some unitary symmetry as in $\mathcal{U}$, see
\eqref{eq:DefU}. In mathematical literature, typically, both
cases are referred to as ~ {\em unitary}.
\vskip 0.2cm

For a given $~p\times m$-valued rational function $F(z)$, let
$F^{\#}(z)$ be the \mbox{$m\times p$-valued} ~{\em conjugate}~ 
rational function, i.e.
\[
F^{\#}(z):=\left(F\left(\frac{1}{z^*}\right)\right)^*.
\]
Note that on the unit circle one has that,
\[
{F^{\#}(z)}_{|_{z\in\mathbb{T}}}=
\left({F(z)}_{|_{z\in\mathbb{T}}}\right)^*.
\]
It is well known,  see e.g. \cite[Eq. (3.1)]{AlGo1}
\cite[Eq. (1.9)]{Oli} that
for rational functions condition \eqref{eq:DefU}
is equivalent to the following, i.e.
\[
\left\{
\begin{matrix}F^{\#}(z)F(z)\equiv I_m&p\geq m&{\rm isometry}\\~\\
F(z)F^{\#}(z)\equiv I_p&m\geq p&{\rm coisometry}.\end{matrix}\right.
\quad\quad\quad\forall z\in\C.
\]
The interest in the class $\mathcal{U}$ is from various aspects, see
e.g. \cite{AlGo1}, \cite{AlGo2}, \cite{AlRak1}, \cite{AlRak2}
\cite{BJ1}, \cite{dBR}, \cite{GVKDM}, \cite{HOP}, \cite{MWBCR}, \cite{Oli},
\cite{PHO}, \cite{Pot}, \cite{SDeLID}, \cite{TAR}, \cite{Valli}.
\vskip 0.2cm

Clearly, whenever $F(z)$ is in $\mathcal{U}$ it must be analytic on
$\mathbb{T}$. There are (at least) two common special cases:
\vskip 0.2cm

(i) If $F(z)$ is analytic outside the closed unit disk
(=Schur stable), then in engineering terminology it is
called ~{\em lossless}\begin{footnote}{Passive electrical
circuits are either dissipative or lossless.}\end{footnote},
see e.g. \cite{GVKDM}, \cite[Section 14.2]{Va}
or ~{\em all-pass}\begin{footnote}{For example, in studying
classical filters a ``high-pass" could be viewed as an
``all-pass" minus a ``low-pass".}\end{footnote}.
\vskip 0.2cm

(ii) If for $p\geq m~$ ($m\geq p$) the matrix
$I_m-(F(z))^*F(z)$~ $~\left(I_p-F(z)(F(z))^*\right)$
is positive semi-definite, within the unit disk, $1\geq|z|$,
then $F(z)$ is anti Schur stable\begin{footnote}{In control engineering
circles a Schur stable functions in $\mathcal{U}$ is called
``inner", see e.g. \cite[Subsection 21.5.1]{ZDG}, while in mathematical
analysis the same term is attributed to the anti Schur stable case,
see e.g. \cite[Section 4]{AlRak2}.}\end{footnote}, i.e.
its conjugate
$F^{\#}(z)$ is Schur stable
\vskip 0.2cm

The interest in 
rational functions within $~{\mathcal U}~$
is vast, see e.g. the books \cite{BJ2}, \cite{J3},
\cite[Section 7.3]{Ma}, \cite[Section 5.2]{SN},
\cite[Section 6.5]{Va} and the papers \cite{AJL1}, \cite{AJL2},
\cite{AJL3}, \cite{AJLM}, \cite{BJ1}, \cite{CZZM4}, \cite{GNS},
\cite{J2}, \cite{J5}, \cite{OTHN}, \cite{RMW}, \cite{TV},
\cite{VHEK} and \cite{YaZh08}.
\vskip 0.2cm

This work is aimed at three different communities: mathematicians
interested in classical analysis, signal processing engineers and
system and control engineers. Thus adopting the terminology
familiar to one audience, may intimidate or even alienate the
other. For example as we already mentioned, rational functions
which are improper or have poles at infinity or non-causal, are
virtually the same entity seen by a different community.
Similarly, what is known to engineers as McMillan degree also
arises in geometry of loop groups as an index. 
\vskip 0.2cm

Books like \cite{BN}, \cite{BJ2}, \cite{SN}, and the theses
\cite{Icart}, \cite{Oli} have made an effort to be at least
``bi-lingual". Lack of space prevents us from providing even a
concise dictionary of relevant terms. Instead, we try to
employ only basic concepts or indicate for references
providing for the necessary background.
\vskip 0.2cm

The differences between scientific communities go beyond terminology.
Closely related problems are formulated not in the same framework.
For example, in many of the engineering references in \eqref{eq:rational}
$F(z)$ is assumed to be analytic outside the open unit disk (=Schur
stable), i.e. $k=0$ and the spectral radius of $A$ is less than one. In
other references $F(z)$ is 
a genuine matrix valued polynomial, i.e. in
\eqref{eq:rational} $B$ or $C$ vanish or in \eqref{eq:poly} $q\geq N$.
We here try to provide a simple, yet full, picture.
\vskip 0.2cm

This work is organized as follows.  
\vskip 0.2cm

In Section \ref{sec:ReqSq} we show that a square rational
$F(z)$ in $~\mathcal{U}~$ can always be truncated (by
eliminating rows or columns) to a rectangular function in
$~\mathcal{U}$. Conversely, a rectangular rational function
in $~\mathcal{U}$, can always be embedded (by adding rows
or columns) in a square function in $~\mathcal{U}$.
\vskip 0.2cm

On the one hand, in the special case where $F(z)$
is analytic outside the open unit disk, this result is well
known. On the other hand if $~\mathcal{U}~$ is substituted
by ~{\em indefinite}~ inner product, this result is not always
true (see discussion below). This suggests that our result
is not trivial.
\vskip 0.2cm

In passing, we explore the controllability and observability
gramians associated with rectangular Schur stable (co)-isometries
on the unit circle.
\vskip 0.2cm

In Section \ref{sec:potapov} we combine
the classical Blaschke-Potapov product formula along
with the main result of the preceding section, to introduce
a characterization of rectangular (co)-isometries
on the unit circle, with poles everywhere (including
infinity) excluding the unit circle.
\vskip 0.2cm

In Section \ref{sec:parametrization} we then exploit the
above characterization to introduce in a compact, convex,
easy-to-use, description of all rational functions in $\mathcal{U}$
parametrized by their McMillan degree and dimensions. Again,
the poles may be everywhere (including infinity) excluding
the unit circle. It is straightforward to restrict this
parametrization to Schur stable functions.
\vskip 0.2cm

This is in particular convenient if one wishes to:\\
(i) Design through
optimization, a rational function (co)-isometric on the unit
circle, see e.g. \cite{CZZM4},
\cite{HTN}, \cite{RMW}, \cite{TV} and \cite{VHEK}.\\
(ii) Iteratively apply para-unitary similarity, see e.g.
\cite[Section 3.3]{Icart}, \cite{MWBCR}, \cite{SDeLID}.
In signal processing literature, this is associated with
with {\em channel equalization}~ and in communications
literature with ~{\em decorrelation of signals}~ or
(iii) Iteratively apply Q-R factorization in the
framework of communications, see e.g. \cite{CB1}
\cite{CB2}.
\vskip 0.2cm

In Section \ref{sec:MFD} we resort to the Matrix Fraction
Description (MFD) of the $p\times m$-valued rational
function $F(z)$, i.e.
\[
F(z)=\left\{
\begin{smallmatrix}
N(z)\left(\Delta(z)\right)^{-1}&&N(z)~~~p\times m-{\rm
valued~polynomial},~~~\Delta(z)~~~m\times m-{\rm valued~polynomial}
&&p\geq m,
\\~\\
\left(\tilde{\Delta}(z)\right)^{-1}\tilde{N}(z)
&&\tilde{N}(z)~~~p\times m-{\rm valued~polynomial},~~~
\tilde{\Delta}(z)~~~p\times p-{\rm valued~polynomial}&&m\geq p.
\end{smallmatrix}\right.
\]
See  e.g. \cite[Chapter 6]{Ka}, \cite[Section 13.3]{Va} or
\cite[Chapter 4]{Vid}. In Theorem \ref{Tm:MFD} we introduce an,
MFD based, easy-to-check characterization of $F(z)$ in
$\mathcal{U}$. Note that this test does not require any
minimality of this representation.
\vskip 0.2cm

In \cite{AJL3} we focus on the subclass rational functions:
In mathematical terms $F(z)$ are $p\times m$-valued polynomials
with powers of possibly mixed signs, i.e. where in
\eqref{eq:rational} the matrix $A$ is nilpotent (i.e. $A^l$
vanishes for some natural $l$). In engineering ``dialect"
these are (not necessarily causal) ~{\em Finite Impulse
Response}~ functions.
We there present three characterizations of those functions
within $\mathcal{U}$.
Here, (in Theorem \ref{Tm:UnitaryPoly} below)
we use Theorem \ref{Tm:MFD} to offer an alternative proof
of one of the main results in \cite{AJL3}.

%
%
%
%
%
%
\section{rectangular vs. square para-unitary
rational functions}
\setcounter{equation}{0}
\label{sec:ReqSq}

In this section we show that in the framework of (co)-isometric rational
functions, the rectangular case is essentially equivalent (in a
rigorous sense, see Theorem \ref{Tm:SqRectU}) to the square case.
\vskip 0.2cm

We do it in two stages. First the easier Schur stable case and then
extend it to rational functions with poles anywhere in the complex
plane (including zero and infinity) but the unit circle.

\subsection{minimal state-space 
realization of Schur stable systems}

This subsection provides known background material used
for the sequel.
\vskip 0.2cm

Recall that if a $p\times m$-valued rational function $F(z)$ is so that
\[
\exists\lim\limits_{z\rightarrow\infty}F(z)
\]
i.e. in \eqref{eq:rational} $k=0$, it is bounded at
infinity\begin{footnote}{in engineering it is colloquially
called ~{\em proper}. Note also that $F(z)$ is referred to
as ~{\em causal}.  This is since that when $F(z)$ is viewed
as the (two-sided) $Z$-transform of a discrete-time sequence
$\Phi(t)$ ($t$ integral variable), then $~\Phi(t)\equiv 0~$
for all $0>t$.}\end{footnote}, then it admits a state space
realization
\begin{equation}\label{eq:RealizFr}
F(z)=C(zI_n-A)^{-1}B+D.
\end{equation}
Sometimes it is convenient to present $~F(z)$ in \eqref{eq:RealizFr}
by its $(n+p)\times(n+m)$ realization matrix $~R$, i.e.
\begin{equation}\label{eq:R}
R:={\footnotesize
\left(\begin{array}{c|c}A&B \\ \hline C&D
\end{array}\right)}.
\end{equation}
A realization is called ~{\em minimal}~ if 
$n$, the dimension of $A$, is the smallest possible.
\vskip 0.2cm

Assuming that $F(z)$ in \eqref{eq:RealizFr} is analytic outside
the open unit disk, in Theorem \ref{Tm:RealizationU} below we
present a characterization, through the corresponding realization
matrix $R$ in \eqref{eq:R}, of Schur stable rectangular
rational functions in $\mathcal{U}$.
\vskip 0.2cm

We here mention some of the existing variants of this result:
The basic case is where $R$ in \eqref{eq:R} is square and the
associated
inner-product is definite. An extension to indefinite inner
product framework appeared in \cite[Theorem 3.1]{AlGo1},
\cite[Theorem 2.1]{AlGo2} and \cite[Lemma 2 \& Theorem 3]{GVKDM}.
In \cite[Theorem 4.5]{AlRak2}, the study was further
generalized to the rectangular case, i.e.
\mbox{$F^*(z)J_pF(z)=J_m$} with $J_p, J_m$ signature matrices,
i.e. diagonal matrices satisfying $J_p^2=I_p$ and $J_m^2=I_m$,
see \cite[Theorem 3.1]{AlRak2}.
\vskip 0.2cm

However, the result in \cite{AlRak2} requires the introduction
of a condition on the ~{\em defect}~ of $F(z)$, for definition
see \cite{For}, \cite[p. 460]{Ka} and
for detailed discussion in the context of rectangular isometries see
\cite[Section 2]{AlRak1}, \cite[Section 2]{AlRak2}. 
\vskip 0.2cm

Restricting the discussion to the Schur stable case (spectrum
within the open unit disk) enabled one to prove the above
result by resorting to a more modest tool from Matrix Theory.
\vskip 0.2cm

\begin{Tm}\label{Tm:RealizationU}
Let $F(z)$ be a $~p\times m$-valued rational function with
poles within the open unit disk (Schur stable).
\vskip 0.2cm

I. Assume that $p\geq m$.
\vskip 0.2cm

(i) $F(z)$ is in $\mathcal{U}$
(=lossless) if and only if,
it admits $(p+n)\times(m+n)$ minimal realization matrix \eqref{eq:R}
\[
R:={\footnotesize
\left(\begin{array}{c|c}A&B \\ \hline C&D
\end{array}\right)}.
\]
satisfying
\begin{equation}\label{eq:SteinRiso}
R^*\cdot{\rm diag}\{I_n\quad I_p\}\cdot{R}={\rm diag}\{I_n\quad I_m\}.
\end{equation}
(ii) If \eqref{eq:SteinRiso} holds, one can always find
$\tilde{B}\in\C^{n\times(p-m)}$ and
$\tilde{D}\in\C^{p\times(p-m)}$ so that the 
$(n+p)\times(n+p)$ 
augmented matrix
\begin{equation}\label{eq:Rnp}
R_{n+p}:={\footnotesize
\left(\begin{array}{c|cc}A&B&\tilde{B}\\ \hline C&D&\tilde{D}
\end{array}\right)},
\end{equation}
is unitary, i.e.
\begin{equation}\label{eq:SteinRnp}
R_{n+p}^*R_{n+p}=I_{n+p}=R_{n+p}R_{n+p}^*~.
\end{equation}
(iii) If \eqref{eq:SteinRnp} holds, one can always find, a constant
isometry $~U_{\rm iso}$ so that
\begin{equation}\label{eq:StaSpaRecSqIso}
{\footnotesize\left(\begin{array}{c|c}A&B\\ \hline C&D
\end{array}\right)}=R=R_{n+p}\cdot
{\footnotesize\left(\begin{array}{c|c}I_n&0_{n\times m}
\\ \hline 0_{p\times n}&U_{\rm iso}
\end{array}\right)}\quad\quad\quad\begin{smallmatrix}
U_{\rm iso}\in\C^{p\times m}
\\~\\U_{\rm iso}^*U_{\rm iso}=I_m
\end{smallmatrix}
\end{equation}
\vskip 0.2cm

II. Assume that $~m\geq p$.

(i) $F(z)$ is in $\mathcal{U}$
(=lossless) if and only if,
it admits $(p+n)\times(m+n)$ minimal realization matrix \eqref{eq:R}
\[
R:={\footnotesize
\left(\begin{array}{c|c}A&B \\ \hline C&D
\end{array}\right)}.
\]
satisfying
\begin{equation}\label{eq:SteinRcoIso}
R\cdot{\rm diag}\{I_n\quad I_m\}\cdot{R^*}={\rm diag}\{I_n\quad I_p\}.
\end{equation}
(ii) If \eqref{eq:SteinRcoIso}  holds, one can always find
$\tilde{C}\in\C^{(m-p)\times n}$ and
$\tilde{D}\in\C^{(m-p)\times m}$ so that the 
$(n+m)\times(n+m)$ augmented
matrix 
\begin{equation}\label{eq:Rnm}
R_{n+m}:={\footnotesize
\left(\begin{array}{c|c}A&B\\ \hline C&D\\ \tilde{C}&\tilde{D}
\end{array}\right)},
\end{equation}
is unitary, i.e.
\begin{equation}\label{eq:SteinRnm}
R_{n+m}^*R_{n+m}=I_{n+m}=R_{n+m}R_{n+m}^*~.
\end{equation}
(iii) If \eqref{eq:SteinRnm} holds, one can always find, a constant
coiometry $~U_{\rm coiso}$ so that
\begin{equation}\label{eq:StaSpaRecSqCoiso}
{\footnotesize\left(\begin{array}{c|c}
A&B\\ \hline C&D\end{array}\right)}
=R={\footnotesize
\left(\begin{array}{c|c}I_n&0_{n\times m}
\\ \hline 0_{p\times n}&U_{\rm coiso}
\end{array}\right)}\cdot{R_{n+m}}\quad\quad\quad
\begin{smallmatrix}
U_{\rm cosio}\in\C^{p\times m}
\\~\\
U_{\rm coiso}U_{\rm coiso}^*=I_p
\end{smallmatrix}
\end{equation}
\end{Tm}
\vskip 0.2cm

{\bf Proof}\quad Assume $p\geq m$
\vskip 0.2cm

Part (i) is an adaption of \cite[Theorem 14.5.1]{Va}.
\vskip 0.2cm

Part (ii) appears in \cite[Lemma 21.21]{ZDG}.
\vskip 0.2cm

Part (iii) follows from the fact that multiplying from
the right a $(n+p)\times(n+p)$ unitary, by a $(n+p)\times(n+m)$
isometry yields another $(n+p)\times(n+m)$ isometry.
\vskip 0.2cm

As the case $m\geq p$ is analogous, its proof is omitted.
\qed
\vskip 0.2cm

As already mentioned, the Schur stable case addressed in Theorem
\ref{Tm:RealizationU},  will be extended to rational functions
with poles anywhere in $\{\C\cup\infty\}\smallsetminus{\mathbb T}$,
in Theorem \ref{Tm:SqRectU} in the next subsection.
\vskip 0.2cm

Still in the Schur stable framework (the spectrum of $A$, the
upper left block of $R$ in \eqref{eq:R} is within the open unit
disk), we now recall the notion of Controllability and Observability
Gramians
(for the continuous-time case see e.g.
\cite[Subsections 9.2.1, 9.2.2]{Ka},
\cite[Sections 3.8, 15.1]{ZDG}
%
We shall denote by
$W_{\rm cont}$, $W_{\rm obs}$, the
$n\times n$
Controllability and
Observability Gramians, respectively, obtained from
the solution to the corresponding Stein equations
\begin{equation}\label{eq:gramian}
W_{\rm cont}-AW_{\rm cont}A^*=BB^*
\quad\quad\quad\quad
W_{\rm obs}-A^*W_{\rm obs}A=C^*C.
\end{equation}
The following, is essentially known, for completeness
a proof is provided.

\begin{Pn}\label{CorollaryGramians}
Let $F(z)$ be a $p\times m$-valued rational function whose poles
are within the open unit disk and denote by
$W_{\rm cont}$, $W_{\rm obs}$ the associated controllability and
observability gramians, respectively.
\vskip 0.2cm

Assume that $F(z)$ is in $\mathcal{U}$.
\vskip 0.2cm

I. If $p\geq m$, $F(z)$ admits a state space realization $R$ in
\eqref{eq:SteinRiso} so that
\[
(I_n-W_{\rm cont})\quad{\rm positive~~semidefinite}
\quad\quad\quad\quad\quad W_{\rm obs}=I_n~.
\]
II. If $m\geq p$, $F(z)$ admits a state space realization $R$ in
\eqref{eq:SteinRcoIso} so that
\[
W_{\rm cont}=I_n\quad\quad\quad\quad\quad
(I_n-W_{\rm obs})\quad{\rm positive~~semidefinite}.
\]
III. If $~p=m$, $F(z)$ admits a state space realization $R$
in \eqref{eq:SteinRiso}, \eqref{eq:SteinRcoIso}
so that
\[
W_{\rm cont}=I_n\quad\quad\quad\quad\quad W_{\rm obs}=I_n~.
\]
\end{Pn}

{\bf Proof}\quad
Indeed, assume $p\geq m$. From the upper left block of 
\eqref{eq:SteinRiso}, it follows that $W_{\rm obs}=I_n$.
Consider now \eqref{eq:Rnp}. The upper left block of
the equation 
\mbox{$R_{n+p}R_{n+p}^*={\rm diag}\{I_n\quad I_p\}$}
reads
\[
I_n-AA^*=BB^*+\tilde{B}\tilde{B}^*.
\]
Now, from \eqref{eq:gramian} we have that
\[
W_{\rm cont}-AW_{\rm cont}A^*=BB^*.
\]
Subtraction of the two equations yields,
\[
(I_n-W_{\rm cont})-A(I_n-W_{\rm cont})A^*=\tilde{B}\tilde{B}^*,
\]
so the first part of the claim is established.
\vskip 0.2cm

As the proof the second part is analogous, it is omitted. The
third part follows from the first two.
\qed
\vskip 0.2cm

We conclude this subsection with a couple of brief comments.
\vskip 0.2cm

(a)\quad Part III of Proposition \ref{CorollaryGramians} is classical,
see e.g. \cite[Section 3]{AlGo1} \cite[Corollary 3]{GVKDM}
and later in \cite[Proposition 1.2.1]{Oli}.
\vskip 0.2cm

(b)\quad The technique employed in \eqref{eq:Rnp} and \eqref{eq:Rnm} in
the proof, is commonly used in system theory for the Hankel norm
approximation and is known as ~{\em all-pass embedding}.

\subsection{Rectangular para-unitary rational functions}
Theorem \ref{Tm:SqRectU}, our first main result, establishes
a close connection between square and rectangular rational
functions in $\mathcal{U}$, with poles at
$\{\C\cup\infty\}\smallsetminus{\mathbb T}$.
\vskip 0.2cm

\begin{Tm}\label{Tm:SqRectU}
Let $F(z)$ be a $p\times m$-valued rational function.

I. Assume that $p\geq m$.
$F(z)$ is in $\mathcal{U}$ if and only if, there exists
in $\mathcal{U}$, a $p\times p$-valued rational function $F_p(z)$, 
so that
\[
F(z)=F_p(z)U_{\rm iso}\quad\quad\quad
\begin{smallmatrix}
U_{\rm iso}\in\C^{p\times m}\\~\\U_{\rm iso}^*U_{\rm iso}=I_m
\end{smallmatrix}
\]
II. Assume that $m\geq p$. 
$F(z)$ is in $\mathcal{U}$ if and only if, there exists 
in $\mathcal{U}$ a $m\times m$-valued rational function $F_m(z)$, 
so that
\[
F(z)=U_{\rm coiso}F_m(z)
\quad\quad\quad\begin{smallmatrix}
U_{\rm coiso}\in\C^{p\times m}
\\~\\U_{\rm coiso}U_{\rm coiso}^*=I_p~.
\end{smallmatrix}
\]
\end{Tm}
\vskip 0.2cm

The proof is relegated further down this subsection.
\vskip 0.2cm

It should be pointed be pointed out that in
\cite[Proposition 2.1]{AlRak1} a similar result is formulated for
the case where on the imaginary axis (instead of the unit circle)
\[
(F(z))^*J_pF(z)=J_m
\]
with $J_m, J_p$ signature matrices, i.e. diagonals satisfying
$J_m^2=I_m$, $J_p^2=I_p$.
\vskip 0.2cm

As already mentioned above, restricting the discussion here to
$J_m=I_m$, $J_p=I_p$ enables us to prove the result through
basic matrix theory tools and to avoid the introduction of the
subtle notion of ~{\em defect} of $F(z)$.
\vskip 0.2cm

In the sequel we shall use the fact that
the scalar rational function (known as a Blaschke-Potapov factor)
\[
\phi(z)=\frac{1-{\alpha}^*z}{z-\alpha}
\quad\quad\quad\quad
\alpha\in\{\infty\cup\C\}\smallsetminus\mathbb{T},
\]
is well defined $\left({\phi(z)}_{|_{\alpha=\infty}}=z\right)$ and
satisfies,
\[
|{\phi(z)}|=1
\quad\quad\quad\quad\forall z\in\mathbb{T}.
\]
We start with an illustrative example.

\begin{Ex}\label{Ex:SqRect1}
{\rm
In part II of Theorem \ref{Tm:SqRectU} take $m=2$,
\begin{equation}\label{eq:ExFm}
F_m(z):={\scriptstyle\frac{1}{\sqrt{2}}}\left(
\begin{smallmatrix}\phi(z)&\psi(z)\\~\\
-\left(\psi(z)\right)^{\#}&\left(\phi(z)\right)^{\#}
\end{smallmatrix}\right),
\end{equation}
where $\phi(z), \psi(z)$ are scalar rational functions.
Then
\[
\left(F_m(z)\right)^{\#}F_m(z)=
{\scriptstyle\frac{1}{2}}\left(
\left(\phi(z)\right)^{\#}\phi(z)
+
\left(\psi(z)\right)^{\#}\psi(z)
\right)I_2~.
\]
Construct from $F_m(z)$ in \eqref{eq:ExFm},
the following $1\times 2$-valued rational function
\begin{equation}\label{eq:ExF1}
F(z)=U_{\rm coiso}F_m(z)
\quad\quad\quad U_{\rm coiso}=(1\quad 0)
\end{equation}
i.e.
\[
F(z)={\scriptstyle\frac{1}{\sqrt{2}}}\left(
\begin{smallmatrix}\phi(z)&&\psi(z)\end{smallmatrix}\right).
\]
Now, $F(z)$ in \eqref{eq:ExF1} is in $\mathcal{U}$, if and only if
$F_m(z)$ in \eqref{eq:ExFm} is in $\mathcal{U}$.
\vskip 0.2cm

This in turn is equivalent to having
$~\phi(z)$, $\psi(z)$ of the form
\[
\phi(z)=\prod\limits_{j=1}^{\overline{j}}
\frac{1-{\alpha_j}^*z}{z-\alpha_j}
\quad\quad
\psi(z)=\prod\limits_{k=1}^{\overline{k}}
\frac{1-{\beta_k}^*z}{z-\beta_k}
\quad\quad\quad\quad\begin{smallmatrix}
\overline{j}, \overline{k}\quad{\rm non-negative~integers}\\~\\
\alpha_j,~\beta_k\in\{\infty\cup\C\}\smallsetminus\mathbb{T}.
\end{smallmatrix}
\]
(Recall, $\prod\limits_1^0:=1$)
}
\qed
\end{Ex}
\vskip 0.2cm

To prove Theorem \ref{Tm:SqRectU}
we resort to the following.

\begin{La}\label{La:poles}
Let $F(z)$ be a $~p\times m$-valued rational function
with poles at \mbox{$\{\infty\cup\C\}\smallsetminus\mathbb{T}$}.
\vskip 0.2cm

I. Assume $p\geq m$
\vskip 0.2cm

One can always find a $m\times m$-valued 
function $U_m(z)$ in $\mathcal{U}$, so that the poles of $F_o(z)$, i.e.
\begin{equation}\label{eq:FU}
F_o(z):=F(z)U_m(z)
\end{equation}
are all in the open unit disk (Schur stable).

Moreover, $F(z)$ is in $\mathcal{U}$, if and only if, $F_o(z)$ is
in $\mathcal{U}$.
\vskip 0.2cm

II. Assume $m\geq p$
\vskip 0.2cm

One can always find a $p\times p$-valued
function $U_p(z)$ in $\mathcal{U}$, so that the poles of $F_o(z)$, i.e.
\[
F_o(z):=U_p(z)F(z)
\]
are all in the open unit disk (Schur stable).

Moreover, $F(z)$ is in $\mathcal{U}$, if and only if, $F_o(z)$ is
in $\mathcal{U}$.
\end{La}
\vskip 0.2cm

{\bf Proof :}\quad
I. Assume $p\geq m$
\vskip 0.2cm

Clearly, for an arbitrary $m\times m$-valued
$U_m(z)$ in $\mathcal{U}$, one has that in \eqref{eq:FU}
$F(z)$ is in $\mathcal{U}$, if and only if, $F_o(z)$ is in
$\mathcal{U}$.
\vskip 0.2cm

Without loss of generality, we shall order the poles of $F(z)$
(including multiplicities) $\alpha_1~,~\ldots~,~
\alpha_t,~\alpha_{t+1},~\ldots~,~\alpha_l$ as
\[
\infty\geq |\alpha_1|\geq\ldots\geq|\alpha_t|>1>
|\alpha_{t+1}|\geq\ldots\geq|\alpha_l|\geq 0.
\]
Take now in \eqref{eq:FU}
\[
U_m(z):=
\prod\limits_{j=1}^t\frac{z-\alpha_j}{1-\alpha_j^*z}I_m~.
\]
It is easy to verify that the poles of $F_o(z)$ in \eqref{eq:FU}
are at
\[
\frac{1}{{\alpha_1}^*}~,~\ldots~,~\frac{1}{{\alpha_t}^*}~,~\alpha_{t+1}~,~
\ldots~,~\alpha_l
\]
and in particular they are all in the open unit disk.
\vskip 0.2cm

The proof of the case $m\geq p$ is analogous and thus omitted.
\qed
\vskip 0.2cm

There are numerous ways to construct $U_m(z)$ in \eqref{eq:FU} (or
$U_p(z)$). The choice in the above proof was solely to simplify
the presentation. It is by no means ``good" in other senses.
\vskip 0.2cm

We can now establish the main result of this section.
\vskip 0.2cm

{\bf Proof of Theorem \ref{Tm:SqRectU}}\quad
If $F(z)$ is Schur stable (poles within the open unit disk), the
claim is established by using $U_{\rm iso}$, $U_{\rm coiso}$
from \eqref{eq:StaSpaRecSqIso},
\eqref{eq:StaSpaRecSqCoiso}, respectively.
\vskip 0.2cm

If the poles of $F(z)$ are anywhere in 
\mbox{$\{\infty\cup\C\}\smallsetminus\mathbb{T}$},
by employing Lemma \ref{La:poles} one may obtain a
Schur stable $F_o(z)$. Now, by the first part, the claim is
established.
\qed
\vskip 0.2cm

The following example illustrates some of the results
of this section.

\begin{Ex}\label{Ex:SqRect2}
{\rm
From Example \ref{Ex:SqRect1} we here consider the 
$1\times 2$-valued $F(z)$ see \eqref{eq:ExF1} and the
$2\times 2$-valued $F_m(z)$ satisfying 
\[
F(z)=(1\quad 0)F_m(z).
\]
For simplicity take in \eqref{eq:ExF1} 
$~\overline{j}=1$, $\overline{k}=0$ so that $F(z)$ and $F_m(z)$ are
of the form
\begin{equation}\label{eq:ExF2}
\begin{matrix}
F(z)&=&
{\scriptstyle\frac{1}{\sqrt{2}}}\left(
\begin{smallmatrix}
\frac{1-{\alpha}^*z}{z-\alpha}&1\end{smallmatrix}\right)\\~\\
F_m(z)&=&
{\scriptstyle\frac{1}{\sqrt{2}}}\left(
\begin{smallmatrix}
\frac{1-{\alpha}^*z}{z-\alpha}&1\\
-1&\frac{z-\alpha}{1-{\alpha}^*z}\end{smallmatrix}
\right)
\end{matrix}
\quad\quad\quad{\scriptstyle\alpha\in\{\infty\cup\C\}\smallsetminus\mathbb{T}}.
\end{equation}
Now, whenever $\alpha$ is restricted to be finite, $F(z)$ in \eqref{eq:ExF1}
admits a (minimal) state space realization of the form \eqref{eq:R} with,
\[
R={\footnotesize\left(\begin{array}{c|cr}
\alpha&\frac{1-|\alpha|^2}{\sqrt{2}}&0\\
\hline
1&-\frac{{\alpha}^*}{\sqrt{2}}&\frac{1}{\sqrt{2}}\end{array}\right)}
\quad\quad\quad{\scriptstyle\alpha\in\{\C\smallsetminus\mathbb{T}\}}.
\]
Furthermore, in accordance to part II of Theorem
\ref{Tm:RealizationU}, it is only when $F(z)$ in \eqref{eq:ExF2} is
lossless (i.e. $1>|\alpha|)$, that it admits an equivalent minimal
realization,
\begin{equation}\label{eq:ExHatR}
\hat{R}={\footnotesize\left(\begin{array}{c|cr}
\alpha&\sqrt{1-|\alpha|^2}&0\\
\hline
\frac{\sqrt{1-|\alpha|^2}}{\sqrt{2}}
&-\frac{{\alpha}^*}{\sqrt{2}}&\frac{1}{\sqrt{2}}
\end{array}\right)},\quad\quad\quad{\scriptstyle 1>|\alpha|},
\end{equation}
satisfying,
\[
\hat{R}\cdot{\rm diag}\{1, I_2\}\cdot\hat{R}^*
={\rm diag}\{1,~1\}.
\]
In fact, following part II of Proposition
\ref{CorollaryGramians}, here the observability gramian is
\mbox{$W_{\rm obs}=\frac{1}{2}$}.
\vskip 0.2cm

Moreover, following \eqref{eq:Rnm}, $\hat{R}$ in \eqref{eq:ExHatR}
may be extended to (here $n=1$, $m=2$),
\[
R_{n+m}={\footnotesize\left(\begin{array}{c|cr}
\alpha&\sqrt{1-|\alpha|^2}&~~~0\\ \hline
\frac{\sqrt{1-|\alpha|^2}}{\sqrt{2}}&-\frac{{\alpha}^*}{\sqrt{2}}
&~~\frac{1}{\sqrt{2}}\\
\frac{\sqrt{1-|\alpha|^2}}{\sqrt{2}}&-\frac{{\alpha}^*}{\sqrt{2}}&
-\frac{1}{\sqrt{2}}
\end{array}\right)},
\]
satisfying
\[
R_{n+m}R_{n+m}^*=I_{m+n}=R_{n+m}^*R_{n+m}.
\]
}
\qed
\end{Ex}

\section{
A characterization through the Blaschke-Potapov product}
\setcounter{equation}{0}
\label{sec:potapov}

We first recall Potapov's classical characterization of the set
of rational functions in $\mathcal{U}$.
Here is a brief perspective. The Fundamental Theorem, see
\cite[p. 133]{Pot}, was formulated in the following framework,
\begin{equation}\label{eq:InDef}
\begin{matrix}
J-F(z)JF^*(z)&{\rm positive~semidefinite}&&1\geq |z|\\~\\
J=F(z)JF^*(z)&&&1=|z|
\end{matrix}\quad\quad\quad\begin{matrix}J\quad{\rm diagonal}\\~\\
J^2=I.\end{matrix}
\end{equation}
A similar result, independently appeared in \cite[Theorem 17]{dBR}
and yet another independent (and more general) version in
\cite[Theorem 9]{GVKDM}. 
\vskip 0.2cm

A special case of this result where $J=I$, was advertized in the
Signal Processing community in \cite[Section 14.9.1]{Va}, see
also \cite{GNS}. In all these cases it was assumed that $F(z)$ is
analytic outside the open unit disk (Schur stable).
\vskip 0.2cm

In \cite[Theorem 3.11]{AlGo1}, Potapov's Fundamental
Theorem was extended to the case where $F(z)$ is analytic on
the circle only (with poles possibly at infinity as well).
\vskip 0.2cm

We shall denote by $P$ a rank one orthogonal projection, i.e.
\[
P^*=P=P^2
\quad\quad\quad\quad{\rm rank}(P)=1.
\]
Recall that if $P$ is $k\times k$ it can always be written as
\begin{equation}\label{eq:DefProjVec}
P=vv^*\quad\quad v^*v=1\quad\quad v\in\C^k.
\end{equation}
Recall also that a rank $k-1$ orthogonal projection $~Q~$ i.e.
\[
Q^*=Q^2=Q\quad\quad\quad\quad\quad\quad{\rm rank}(Q)=k-1,
\]
can always be written as
\begin{equation}\label{eq:Qprojection}
Q:=I_k-vv^*\quad\quad v^*v=1\quad\quad v\in\C^k
\end{equation}
as in \eqref{eq:DefProjVec}.

\begin{Tm}\label{Tm:U(z)}
\mbox{Let $F(z)$ be a $~p\times m$-valued rational function of McMillan
degree $~d$.}

\vskip 0.2cm
$F(z)$ is in $~{\mathcal U}$, \eqref{eq:DefU}, if and only if it
can be written as
\begin{equation}\label{eq:UnitaryRational1}
\begin{matrix}
p\geq m&F(z)=\left(\prod\limits_{j=1}^d\left(I_p+\left(
\frac{1-{\alpha}^*_jz}{z-\alpha_j}-1\right)v_jv_j^*\right)\right)
U_{\rm iso}\\~\\
m\geq p&F(z)=U_{\rm coiso}\left(\prod\limits_{j=1}^d
\left(I_m+\left(\frac{1-{\alpha}^*_jz}{z-\alpha_j}-1\right)v_jv_j^*
\right)\right)
\end{matrix}\quad
\begin{smallmatrix}
v_j\in\C^p&
\quad v_j^*v_j=1
\\~\\
U_{\rm iso}\in\C^{p\times m}&
U_{\rm iso}^*U_{\rm iso}=I_m
\\~\\
\quad\alpha_j\in\{\infty\cup\C\}\smallsetminus\mathbb{T}&~\\~\\
v_j\in\C^m&v_j^*v_j=1\\~\\
U_{\rm coiso}\in\C^{p\times m}&U_{\rm coiso}U_{\rm coiso}^*=I_p~.
\end{smallmatrix}
\end{equation}
\end{Tm}
\vskip 0.2cm

Recall $\prod\limits_{j=1}^0:=I$
\vskip 0.2cm

{\bf Proof}\quad
Substituting in  \cite[Theorem 3.11]{AlGo1} the special case
$J=I$ (definite inner product), yields the following:

An $~m\times m$-valued rational function $F(z)$, of McMillan
degree $~d$, is in $~{\mathcal U}$, \eqref{eq:DefU}, if and
only if (up to multiplication by a constant $m\times m$
unitary matrix from the left or from the right) it can be
written as
\begin{equation}\label{eq:UnitaryRational}
F(z)=\prod\limits_{j=1}^d
\left(I_m+\left(\frac{1-{\alpha}^*_jz}{z-\alpha_j}
-1\right)v_jv_j^*\right)
\quad\quad\quad
\alpha_j\in\{\infty\cup\C\}\smallsetminus\mathbb{T}.
\end{equation}
Using, \eqref{eq:DefProjVec} and \eqref{eq:Qprojection},
establishes \eqref{eq:UnitaryRational1} for $~m=p$.
\vskip 0.2cm

To obtain the rectangular case, apply Theorem \ref{Tm:SqRectU}.
\qed
\vskip 0.2cm

Three remarks are now in order.
\vskip 0.2cm

{\bf a}\quad 
It is tempting to combine \cite[Theorem 3.11]{AlGo1} along
with the above Theorem \ref{Tm:U(z)}, to formulate a rectangular
version of Blaschke-Potapov product result with poles in
$\{\infty\cup\C\}\smallsetminus\mathbb{T}$ for ~{\em indefinite}~ 
inner product, see \eqref{eq:InDef}. However, this requires some
caution as then, the notion of the ~{\em defect}~ of $F(z)$ needs
to be addressed. For definition see \cite{For}, \cite[p. 460]{Ka} and for detailed
discussion in the context of rectangular isometries see
\cite[Section 2]{AlRak1}.
\vskip 0.2cm

{\bf b}\quad
Theorem \ref{Tm:U(z)} asserts that whenever $F\in\mathcal{U}$
is of McMillan degree $d$, {\em there exist}~ rank one orthogonal
projections $~v_1v_1^*~,~\ldots~,~v_dv_d^*$, satisfying
\eqref{eq:UnitaryRational1}. In general, the McMillan degree of
the product in the right hand side of \eqref{eq:UnitaryRational}
is ~{\em at most}~ $d$. For example,
\[
\left(I+(\phi_1(z)-1)v_1v_1^*\right)
\left(I+(\phi_2(z)-1)v_2v_2^*\right)_{|_{v_1v_1^*=v_2v_2^*}}
=\left(I+(\phi_1(z)\phi_2(z)-1)v_1v_1^*
\right)_{|_{\phi_1(z)\phi_2(z)\equiv 1}}=I,
\]
which is a zero degree rational function.
\vskip 0.2cm

{\bf c}\quad 
Note that products of the form 
\[
v_1v_1^*v_2v_2^*\cdots{v_kv_k^*}=
\left(\prod\limits_{j=1}^{k-1}v_j^*v_{j+1}\right)v_1v_k^*
\quad\quad\quad\quad k\geq 2,
\]
which appear in \eqref{eq:UnitaryRational1}, always produce a rank
one matrix. In the special case where $~v_1v_1^*=~\cdots~=v_dv_d^*$
this is an orthogonal
projection, else it is a strict contraction (which may be Hermitian
when $v_1$ and $v_k$ are linearly dependent).
\vskip 0.2cm

%
%
%

\section{parametrization of all para-unitary rational functions}
\setcounter{equation}{0}
\label{sec:parametrization}

We next exploit the above Theorem \ref{Tm:U(z)} to describe
all rational function in ${\mathcal U}$, parametrized by 
dimensions and the McMillan degree.
\vskip 0.2cm

To this end, we introduce the following matrix theory notation
\begin{equation}\label{eq:DefIso}
\begin{matrix}
{\mathbb U}_{\rm Iso}&:=&\{ U\in\C^{p\times m}&p\geq m~:~U^*U=I_m~\}\\~\\
{\mathbb U}_{\rm Coiso}&:=&\{ U\in\C^{p\times m}&m\geq p~:~UU^*=I_p~\}.
\end{matrix}
\end{equation}
\begin{La}\label{Lemma:Iso}
The set ${\mathbb U}_{\rm Iso}$ in \eqref{eq:DefIso} may be completely
parametrized by
\begin{equation}\label{eq:ParameterIso}
[0,~2\pi)^{m(2p-m)}.
\end{equation}
Similarly, the set ${\mathbb U}_{\rm Coiso}$ in \eqref{eq:DefIso} may
be completely parametrized by
\[
[0,~2\pi)^{p(2m-p)}.
\]
\end{La}

Indeed, due to symmetry, one address
only the case of $p\geq m$. Now, 
the set of all $v\in\C^p$ with $v^*v=1$, i.e. the $\|~\|_2$ unit
sphere in $\in\C^p$ may be identified with with
\[
[0,~2\pi)^{2p-1}.
\]
For example for $p=3$ this $v$ is of the form
\[
v=\left(\begin{smallmatrix}
\cos(\alpha)e^{i\eta}\\ \cos(\beta)\sin(\alpha)e^{i\gamma}\\
\sin(\beta)\sin(\alpha)e^{i\delta}\end{smallmatrix}\right)
\quad\quad\quad\quad\quad
\alpha, \beta, \gamma, \delta, \eta\in[0,~2\pi)
\]
To obtain all $m$-dimensional orthonormal bases of
such vectors, one resorts to \eqref{eq:ParameterIso}, so
the claim is established.
\vskip 0.2cm

A word of caution. Consider for simplicity the case of unitary
matrices where $p=m$ are prescribed. One can ask the two
following questions. 

(i) How many
parameters are required to completely describe the whole set.\\
(ii) How many parameters are required to completely describe all
unitary similarity transformations.
\vskip 0.2cm

The above Lemma addresses the first question.
The following example illustrates the gap between these two.
\vskip 0.2cm

\begin{Ex}
{\rm Consider for simplicity the case of $p=m=2$.
\vskip 0.2cm

Every unitary matrix $~U~$ may be written as
\[
U=\left(\begin{smallmatrix}
e^{i(\gamma-\beta)}\cos(\alpha)&~~~~~~~~~~~~
e^{i\delta}\sin(\alpha)\\~~~~~~~-e^{-i\beta}\sin(\alpha)&~~~
e^{i(\delta-\gamma)}\cos(\alpha)
\end{smallmatrix}\right)\quad\quad
\alpha,\beta,\gamma,\delta\in[0,~2\pi).
\]
Namely, this set may be identified with $[0,~2\pi)^4$.
\vskip 0.2cm

However, if for a given $2\times 2$ matrix $M$, one is interested
in all unitary similarity transformations of the form $U^*MU$,
without loss of generality, one can assume that in the above $U$,
\[
\beta=\gamma=\delta.
\]
Namely, two of the angles are redundant, so 
all $2\times 2$ unitary similarity transformations
may be identified with $[0,~2\pi)^2.$
\vskip 0.2cm

In this case the complex version of the Givens (sometimes named
after Jacobi) rotations is obtained
(for the real version see e.g. \cite[Section 3.4]{GVL},
\cite[Example 2.2.3]{HJ1}
\cite[Section 14.6.1]{Va}). Thus, it is parametrized by two
(and not four) angles. 
\vskip 0.2cm

In the literature these two problems were treated in numerous places
(in some cases, with a slight confusion between them), see e.g.
\cite[Section 3.4]{GVL}, \cite[Propri\'{e}t\'{e} 41]{Icart}, \cite{MWBCR},
\cite[Eq. (19)]{RMW}, \cite[Section 3]{SDeLID}, \cite{TAR} and
\cite[Section 14.6.1]{Va}.
}
\qed
\end{Ex}
\vskip 0.2cm

Theorem \ref{Tm:U(z)} along with Lemma \ref{Lemma:Iso}
enable us to introduce the following easy-to-use description of all
rational functions in $\mathcal{U}$ of prescribed McMillan degree
$d$ and dimensions $p$ and $m$, as real set which is virtually
$d$ copies of real polytopes.

\begin{Ob}\label{Ob:angles}
All $p\times m$-valued rational functions of McMillan degree
$~d~$ in $\mathcal{U}$
may be parametrized by,
\begin{equation}\label{eq:angle}
\begin{matrix}
\left(\{0\}\cup\{\infty\}\cup
\left(\left( (0,~\infty)\smallsetminus\{ 1\}\right)\cdot[0,~2\pi)\right)
\right)^d\cdot[0,~2\pi)^{2d(p-1)+m(2p-m)}
&p\geq m\\~\\
\left(\{0\}\cup\{\infty\}\cup
\left(\left( (0,~\infty)\smallsetminus\{ 1\}\right)\cdot[0,~2\pi)\right)
\right)^d\cdot[0,~2\pi)^{2d(m-1)+p(2m-p)}
&m\geq p.
\end{matrix}
\end{equation}
The Schur stable subset is parametrized by,
\[
\begin{matrix}
\left(\{0\}\cup
\left((0,~1)\cdot[0,~2\pi)\right)
\right)^d\cdot[0,~2\pi)^{2d(p-1)+m(2p-m)}
&p\geq m\\~\\
\left(\{0\}\cup
\left((0,~1)\cdot[0,~2\pi)\right)
\right)^d\cdot[0,~2\pi)^{2d(m-1)+p(2m-p)}
&m\geq p.
\end{matrix}
\]
\end{Ob}

{\bf Proof:}\quad
Assume that $p\geq m$. As in Lemma \ref{Lemma:Iso} the set of all
$v\in\C^p$ with $v^*v=1$,
i.e. the $\|~\|_2$ unit sphere in $\C^p$, may be identified with
\[
[0,~2\pi)^{2p-1}.
\]
As $~v$ and $e^{i\eta}v~$ produce the same $vv^*$, to
parametrize all $p\times p$ rank one orthogonal projections
in \eqref{eq:DefProjVec}, one angle is redundant, so one can use
\[
[0,~2\pi)^{2(p-1)}.
\]
We next address the poles $\alpha_1~,~\ldots~,~\alpha_d$ in
\eqref{eq:UnitaryRational1}. If a pole $\alpha_j$ is in the
complex plane, excluding zero, infinity and the unit circle,
it may be parametrized by the usual polar representation,
\begin{equation}\label{eq:FinitePoles}
\left( (0,~\infty)
\smallsetminus\{ 1\}
\right)\cdot[0,~2\pi).
\end{equation}
Thus, to parametrize a single Blaschke-Potapov factor in
\eqref{eq:UnitaryRational1}, one needs
\[
\left(
\{0\}\cup\{\infty\}\cup\left( (0,~\infty)\smallsetminus\{ 1\}\right)\cdot[0,~2\pi)
\right)\cdot[0,~2\pi)^{2(p-1)}.
\]
Note that this set is nearly a real polytope.
Now taking $d$ copies, yields,
\[
\left(
\{0\}\cup\{\infty\}\cup\left( (0,~\infty)\smallsetminus\{ 1\}\right)\cdot[0,~2\pi)
\right)^d\cdot[0,~2\pi)^{2d(p-1)}.
\]
Along with \eqref{eq:ParameterIso} from Lemma \ref{Lemma:Iso} the
first part of \eqref{eq:angle} is obtained.
\vskip 0.2cm

Due to symmetry, we omit the case $m\geq p$, so the construction
is complete.
%
%
%
%
%
\qed
\vskip 0.2cm

The above parameterization is in particular convenient if one
wishes to design through
optimization, a rational function (co)-isometric on the unit
circle. For example, given a $p\times m$-valued function $G(z)$
which is not necessarily rational, not necessarily (co)-isometric
on the unit circle, and not necessarily Schur stable, find
$F(z)$ its best Schur stable approximation in $\mathcal{U}$ of a 
prescribed McMillan degree $~d$, i.e.
\[
\min\limits_{
\left(\{0\}\cup
\left((0,~1)\cdot[0,~2\pi)\right)
\right)^d\cdot[0,~2\pi)^{2d(p-1)+m(2p-m)}
}
\| F(z)-G(z)\| \quad\quad\quad\quad p\geq m.
\]
For other type optimization problems see e.g. \cite{CZZM4},
\cite{HTN}, \cite{RMW}, \cite{TV} and \cite{VHEK}.

\section{matrix-fraction description}
\setcounter{equation}{0}
\label{sec:MFD}

So far, we confined the discussion to rational functions
$F(z)$ presented in their ~{\em minimal realization}. 
We next relax this restriction.
\vskip 0.2cm

Following e.g. \cite[Chapter 6]{Ka}, \cite[Section 13.3]{Va}
or \cite[Chapter 4]{Vid}, a $~p\times m$-valued rational function
of the form \eqref{eq:rational}, can always be written as
\begin{equation}\label{eq:MFDorig}
F(z)=\left\{
\begin{smallmatrix}
N(z)\left(\Delta(z)\right)^{-1}&=&
\left(N_0+zN_1+~\cdots~+z^{\nu}N_{\nu}\right)
\left(\Delta_0+z\Delta_1+~\cdots~+z^{\delta}\Delta_{\delta}\right)^{-1}
&~&{\rm RMFD}
\\~\\
\left(\tilde{\Delta}(z)\right)^{-1}\tilde{N}(z)
&=&
\left(\tilde{\Delta}_0+z\tilde{\Delta}_1+~\cdots~+
z^{\tilde\delta}\tilde{\Delta}_{\tilde\delta}\right)^{-1}
\left(\tilde{N}_0+z\tilde{N}_1+~\cdots~+z^{\tilde\nu}
\tilde{N}_{\tilde\nu}\right)&~&{\rm LMFD}
\end{smallmatrix}\right.
\end{equation}
where $\Delta(z)$ and $\tilde{\Delta}(z)$ are $~m\times m$-valued and
$~p\times p$-valued polynomials, respectively, each of a full normal
rank, while both $N(z)$ and $\tilde{N}(z)$ are $p\times m$-valued
polynomials.  $N(z)\left(\Delta(z)\right)^{-1}$ is called a ~{\em right
matrix fraction description} (RMFD) of $F(z)$ while
$\left(\tilde{\Delta}(z)\right)^{-1}\tilde{N}(z)$
is a ~{\em left matrix fraction description}~
(LMFD) of $F(z)$.
\vskip 0.2cm

Specifically, $\nu$, $\delta$, $\tilde{\nu}$ and $\tilde{\delta}$ in
\eqref{eq:MFDorig} are non-negative integers. If they are the smallest
possible\begin{footnote}{In principle,
for arbitrary $m\times m$-valued polynomial $R(z)$, another
RMFD is \mbox{$F(z)=N(z)R(z)\left(\Delta(z)R(z)\right)^{-1}$}}
\end{footnote} the matrix
fraction description of $F(z)$ in \eqref{eq:MFDorig} is said to be~
{\em irreducible}~ see e.g. \cite[subsection 6.5]{Ka}. Then,
the polynomials $N(z)$ and $\Delta(z)$ are right coprime or the
polynomials $\tilde{N}(z)$ and $\tilde{\Delta}(z)$ are left coprime,
for details see e.g. \cite[subsection 6.5]{Ka} or \cite[Chapter 4]{Vid}.
\vskip 0.2cm

For a given $F(z)$, finding an {\em irreducible} MFD, may be
challenging.  However, here we look for ~{\em some} MFD.
Specifically, let,
\[
\alpha\geq\max(\nu,~\delta)
\quad\quad\quad\quad
\beta\geq\max(\tilde{\nu},~\tilde{\delta}),
\]
and by formally adding zero matrices to \eqref{eq:MFDorig}, we
shall hereafter use the following MFD, where the numerator
and denominator polynomials have the same power,
\begin{equation}\label{eq:MFD}
F(z)=\left\{
\begin{smallmatrix}
N(z)\left(\Delta(z)\right)^{-1}&=&
\left(N_0+zN_1+~\cdots~+z^{\alpha}N_{\alpha}\right)
\left(\Delta_0+z\Delta_1+~\cdots~+z^{\alpha}\Delta_{\alpha}\right)^{-1}
\\~\\
\left(\tilde{\Delta}(z)\right)^{-1}\tilde{N}(z)
&=&
\left(\tilde{\Delta}_0+z\tilde{\Delta}_1+~\cdots~+
z^{\beta}\tilde{\Delta}_{\beta}\right)^{-1}
\left(\tilde{N}_0+z\tilde{N}_1+~\cdots~+z^{\beta}
\tilde{N}_{\beta}\right)
\end{smallmatrix}\right.
\end{equation}
Recall also that with the polynomials in the RMFD in \eqref{eq:MFD}
one can associate the following ($p(\alpha+1)\times m(\alpha+1)~$
and $~m(\alpha+1)\times m(\alpha+1)$, respectively) Hankel matrices
\begin{equation}\label{eq:HankelRight}
{\mathbf H}_N:=\left(\begin{smallmatrix}
N_0    &N_1   & ~    &N_{\alpha-1}&N_{\alpha}\\
N_1    &  ~   &\Ddots&\Ddots &~ \\
~      &\Ddots&\Ddots&~     &~  \\
N_{\alpha-1}&\Ddots&~     &~     &~\\
N_{\alpha}&~     &~     &~     &\end{smallmatrix}\right)
\quad\quad\quad\quad
{\mathbf H}_{\Delta}:=\left(\begin{smallmatrix}
\Delta_0    &\Delta_1   & ~    &\Delta_{\alpha-1}&\Delta_{\alpha}\\
\Delta_1    &  ~   &\Ddots&\Ddots &~ \\
~      &\Ddots&\Ddots&~     &~  \\
\Delta_{\alpha-1}&\Ddots&~     &~     &~\\
\Delta_{\alpha}&~     &~     &~     &\end{smallmatrix}\right).
\end{equation}
By construction, both ${\mathbf H}_N^*{\mathbf H}_N$ and
${\mathbf H}_{\Delta}^*{\mathbf H}_{\Delta}$ are of
the same dimensions \mbox{$m(\alpha+1)\times m(\alpha+1)$}.
\vskip 0.2cm

Similarly, 
with the polynomials in the LMFD in \eqref{eq:MFD}, one can
associate the following ($p(\beta+1)\times m(\beta+1)~$ and
$~p(\beta+1)\times p(\beta+1)$, respectively) Hankel matrices,
\begin{equation}\label{eq:HankelLeft}
\mathbf{H}_{\tilde{N}}:=\left(\begin{smallmatrix}
\tilde{N}_0    &\tilde{N}_1   & ~    &\tilde{N}_{\beta-1}&\tilde{N}_{\beta}\\
\tilde{N}_1    &  ~   &\Ddots&\Ddots &~ \\
~      &\Ddots&\Ddots&~     &~  \\
\tilde{N}_{\beta-1}&\Ddots&~     &~     &~\\
\tilde{N}_{\beta}&~     &~     &~     &\end{smallmatrix}\right)
\quad\quad\quad\quad
\mathbf{H}_{\tilde\Delta}:=\left(\begin{smallmatrix}
\tilde{\Delta}_0    &\tilde{\Delta}_1   & ~    &
\tilde{\Delta}_{\beta-1}&\tilde{\Delta}_{\beta}\\
\tilde{\Delta}_1    &  ~   &\Ddots&\Ddots &~ \\
~      &\Ddots&\Ddots&~     &~  \\
\tilde{\Delta}_{\beta-1}&\Ddots&~     &~     &~\\
\tilde{\Delta}_{\beta}&~     &~     &~     &\end{smallmatrix}\right).
\end{equation}
By construction, both ${\mathbf H}_{\tilde{N}}{\mathbf H}_{\tilde{N}}^*$
and
${\mathbf H}_{\tilde{\Delta}}{\mathbf H}_{\tilde{\Delta}}^* $ are of
the same dimensions \mbox{$p(\beta+1)\times p(\beta+1)$}.
\vskip 0.2cm

We can now state the main result of this section.

\begin{Tm}\label{Tm:MFD}
Let $F(z)$ be a $p\times m$-valued rational function with a (not
necessarily reducible) Matrix Fraction Description in \eqref{eq:MFD}.
\vskip 0.2cm

(I). For $p\geq m~$ let ${\mathbf H}_N$ and $~{\mathbf H}_{\Delta}$
in \eqref{eq:HankelRight} be the Hankel matrices associated with
$~N(z)$ and $\Delta(z)$, respectively.\quad
%
$F(z)$ is in $\mathcal{U}$, if and only if,
\begin{equation}\label{eq:CondH*H}
\left({\mathbf H}_{\Delta}^*{\mathbf H}_{\Delta}-
{\mathbf H}_N^*{\mathbf H}_N\right)
\left(\begin{smallmatrix}I_m\\ 0_{m\alpha\times m}
\end{smallmatrix}\right)=0_{m(\alpha+1)\times m}~.
\end{equation}
(II). For $m\geq p$ let 
${\mathbf H}_{\tilde{N}}$ and $~{\mathbf H}_{\tilde{\Delta}}$
in \eqref{eq:HankelLeft} be the Hankel matrices associated with
$\tilde{N}(z)$ and $\tilde{\Delta}(z)$, respectively.\quad
%
$F(z)$ is in $\mathcal{U}$, if and only if,
\[
\left({\mathbf H}_{\tilde{\Delta}}{\mathbf H}_{\tilde{\Delta}}^*
-{\mathbf H}_{\tilde{N}}{\mathbf H}_{\tilde{N}}^*\right)
\left(\begin{smallmatrix}I_p\\ 0_{p\beta\times p}
\end{smallmatrix}\right)=0_{p(\beta+1)\times p}~.
\]
\end{Tm}

{\bf Proof :}\quad
Assume that ~$p\geq m$, take the right RMFD of $F(z)$
and consider the following
(where to simplify the presentation we
omit the explicit dependence on the variable $z$)
\[
F^{\#}F=\left(N\Delta^{-1}\right)^{\#}N\Delta^{-1}
=\left(\Delta^{-1}\right)^{\#}N^{\#}N\Delta^{-1}
\]
Now, having $F(z)$ is in $\mathcal{U}$ is equivalent to
\[
I_m=F^{\#}F=\left(\Delta^{-1}\right)^{\#}N^{\#}N\Delta^{-1}.
\]
Multiplying by $\Delta^{\#}$ from the left and $\Delta$ from
the right yields 
\[
\Delta^{\#}\Delta=N^{\#}N.
\]
Substituting now \eqref{eq:MFD} in the above reads
\[
\begin{matrix}
&
\left(\Delta_0+z\Delta_1+~\cdots~+z^{\alpha}\Delta_{\alpha}\right)^{\#}
\left(\Delta_0+z\Delta_1+~\cdots~+z^{\alpha}\Delta_{\alpha}\right)\\
=&\\
&\left(N_0+zN_1+~\cdots~+z^{\alpha}N_{\alpha}\right)^{\#}
\left(N_0+zN_1+~\cdots~+z^{\alpha}N_{\alpha}\right)
\end{matrix}
\]
which is equal to
\[
\begin{matrix}
&
\left(\Delta_0^*+\frac{1}{z}\Delta_1^*+~\cdots~+
\frac{1}{z^{\alpha}}\Delta_{\alpha}^*\right)
\left(\Delta_0+z\Delta_1+~\cdots~+z^{\alpha}\Delta_{\alpha}\right)\\
=&\\
&\left(N_0^*+\frac{1}{z}N_1^*+~\cdots~+
\frac{1}{z^{\alpha}}N_{\alpha}^*\right)
\left(N_0+zN_1+~\cdots~+z^{\alpha}N_{\alpha}\right).
\end{matrix}
\]
Note that in both, the numerator and the denominator, for each
$k\in[1, \alpha]$, the coefficient of $\frac{1}{z^k}$, is the
complex conjugate transpose, $(~)^*$, of the coefficient of $z^k$.
Thus, without loss of generality, one can equate only the
coefficients of $z^k$ for $k\in[0,~\alpha]$. This means that
\[
\mathbf{H}_{\Delta}^*\left(\begin{smallmatrix}
\Delta_o\\ \vdots\\~\\ \Delta_{\alpha}
\end{smallmatrix}\right)=
\mathbf{H}_N^*\left(\begin{smallmatrix}
N_o\\ \vdots\\~\\ N_{\alpha}\end{smallmatrix}\right),
\]
with the Hankel matrices from \eqref{eq:HankelRight}.
This in turn may be equivalently written as
\[
\mathbf{H}_{\Delta}^*\mathbf{H}_{\Delta}
\left(\begin{smallmatrix}I_m\\ 0_{m\delta\times m}\end{smallmatrix}\right)=
\mathbf{H}_N^*\mathbf{H}_N
\left(\begin{smallmatrix}I_m\\ 0_{m\delta\times m}\end{smallmatrix}\right),
\]
so \eqref{eq:CondH*H} is established.
\vskip 0.2cm

Due to symmetry, establishing the case $~m\geq p$, is
analogous and thus omitted.
\qed
\vskip 0.2cm

This work is devoted to $p\times m$-valued {\em rational functions}~
within $\mathcal{U}$. In \cite{AJL3} we focused on the subset of
(possibly Laurent) polynomials (within $\mathcal{U}$) i.e.
\begin{equation}\label{eq:poly}
F(z)=z^q(B_o+zB_1+~\cdots~+z^{\gamma}B_{\gamma})
\quad\quad\gamma~~{\rm natural},~~q~~{\rm integral~~parameter},
\end{equation}
and $B_o, B_1,~\cdots~,~B_{\gamma}$ constant matrices\begin{footnote}
{Strictly speaking, the notation in \cite{AJL3} is slightly
different, but equivalent.}\end{footnote}.
Note that for $-1\geq q$ this is no longer a genuine
polynomial. Although modest is size, there is a vast
literature on this family, see e.g. \cite{AJL3} and references therein.
\vskip 0.2cm

In Theorem \ref{Tm:UnitaryPoly} below we show how to use Hankel
matrices to characterize this subset. In fact, this is a citation
of \cite[theorem 4.1]{AJL3}. However, as the original proof is
somewhat different. Using the above Theorem \ref{Tm:MFD}, we
next establish the same result independently.
\vskip 0.2cm

Here are the details:
Substituting $q=0$ in \eqref{eq:poly} one obtains,
\[
F_0(z):={F(z)}_{|_{q=0}}=B_o+zB_1+~\cdots~+z^{\gamma}B_{\gamma}~.
\]
With $F_0(z)$ one can associate the following \mbox{$p(\gamma+1)\times m(\gamma+1)$}
Hankel matrix,
\begin{equation}\label{eq:Ho}
{\mathbf H}_0:=\left(\begin{smallmatrix}
B_0    &B_1   & ~    &B_{\gamma-1}&B_{\gamma}\\
B_1    &  ~   &\Ddots&\Ddots &~ \\
~      &\Ddots&\Ddots&~     &~  \\
B_{\gamma-1}&\Ddots&~     &~     &~\\
B_{\gamma}&~     &~     &~     &\end{smallmatrix}\right).
\end{equation}

\begin{Tm}\label{Tm:UnitaryPoly}
Let $F(z)$ be a $p\times m$ polynomial in \eqref{eq:poly} and let
${\mathbf H}_0$ be the associated Hankel matrix as in \eqref{eq:Ho}.
\vskip 0.2cm
 
The polynomial $F(z)$ is in $\mathcal{U}$, if and only if,
\begin{equation}\label{eq:H^*H}
\begin{matrix}
\left(I_{m(\gamma+1)}-{\mathbf H}_0^*{\mathbf H}_0\right)
\left(\begin{smallmatrix}
I_m\\0_{m\gamma\times m}\end{smallmatrix}\right)&=
&0_{m(\gamma+1)\times m}&&&p\geq m\\~\\
\left(I_p\quad 0_{p\times p\gamma}\right)
\left(I_{p(\gamma+1)}-{\mathbf H}_0{\mathbf H}_0^*\right)
&=&0_{p\times p(\gamma+1)}&&&m\geq p.
\end{matrix}
\end{equation}
\end{Tm}

{\bf Proof}\quad First, note that if $F(z)$ in \eqref{eq:poly}
is in $\mathcal{U}$ for ~{\em some}~ $q$, it is in $\mathcal{U}$ for ~
{\em all}~ $q$. Thus, without loss of generality, we characterize
$F_0(z)$ in $\mathcal{U}$.
\vskip 0.2cm
 
First, note that as a rational function $F_o(z)$ can be written
as a RMFD in \eqref{eq:MFD} with $\Delta_0=I_m$,
$\Delta_1=~\cdots~=\Delta_{\gamma}=0$
and $N_j=B_j$ for $j=0,~\ldots~,~\gamma$. Thus, using ${\mathbf H}_0$
from \eqref{eq:Ho} here \eqref{eq:HankelRight} takes the form
\[
{\mathbf H}_N={\mathbf H}_0\quad\quad\quad\quad
{\mathbf H}_{\Delta}=\left(\begin{smallmatrix}I_m&&0\\0&&0_{m\gamma\times
m\gamma}\end{smallmatrix}\right).
\]
Thus, for $p\geq m$ using \eqref{eq:CondH*H} one has that,
\[
\begin{smallmatrix}
\left(I_{m(\gamma+1)}-{\mathbf H}_0^*{\mathbf H}_0\right)
\left(\begin{smallmatrix}
I_m\\0_{m\gamma\times m}\end{smallmatrix}\right)&=&
\left(I_{m(\gamma+1)}-{\mathbf H}_N^*{\mathbf H}_N\right)
\left(\begin{smallmatrix}
I_m\\0_{m\gamma\times m}\end{smallmatrix}\right)\\~\\~&=&
\left(I_{m(\gamma+1)}-{\mathbf H}_N^*{\mathbf H}_N+
{\mathbf H}_{\Delta}^*{\mathbf H}_{\Delta}
-{\mathbf H}_{\Delta}^*{\mathbf H}_{\Delta}\right)
\left(\begin{smallmatrix}
I_m\\0_{m\gamma\times m}\end{smallmatrix}\right)\\~\\~&=&
\left(I_{m(\gamma+1)}-{\mathbf H}_{\Delta}^*{\mathbf H}_{\Delta}
\right)\left(\begin{smallmatrix}I_m\\0_{m\gamma\times m}\end{smallmatrix}\right)
+\left({\mathbf H}_{\Delta}^*{\mathbf H}_{\Delta}-
{\mathbf H}_N^*{\mathbf H}_N\right)\left(\begin{smallmatrix}I_m\\0_{m\gamma\times
m}\end{smallmatrix}\right)\\~\\~&=&
\left(I_{m(\gamma+1)}-{\mathbf H}_{\Delta}^*{\mathbf H}_{\Delta}
\right)\left(\begin{smallmatrix}I_m\\0_{m\gamma\times m}
\end{smallmatrix}\right)\\~\\~&=&
\left(\left(\begin{smallmatrix}I_m&&0\\0&&I_{m\gamma}\end{smallmatrix}\right)
-\left(\begin{smallmatrix}I_m&&0\\0&&0_{m\gamma\times m\gamma}
\end{smallmatrix}\right)\right)\left(\begin{smallmatrix}I_m\\0_{m\gamma\times m}
\end{smallmatrix}\right)\\~\\~&=&0_{m(\gamma+1)\times m}
\end{smallmatrix}
\]
Thus, the first part of \eqref{eq:H^*H} is obtained.
\vskip 0.2cm

Due to symmetry, establishing the case $~m\geq p~$ is analogous and thus
omitted.
\qed
\vskip 0.2cm

\bibliographystyle{plain}

\end{document}